\documentclass[12pt]{article}

\usepackage{lipsum}

\usepackage{tikz-cd}
\usepackage{amsmath}
\usepackage{amssymb}
\usepackage{amsthm, thmtools}
\usepackage{mathtools}
\usepackage{mathrsfs}

\usepackage{bbm}
\usepackage{bm}
\usepackage{MnSymbol}   
\usepackage{cases}

\usepackage[width=.80\textwidth,font=footnotesize]{caption}
\usepackage[letterpaper, margin=1in]{geometry}
\usepackage{appendix}
\usepackage{fancyhdr}
\usepackage{caption}
\usepackage{subcaption}

\usepackage[shortlabels]{enumitem}
\usepackage[ddmmyyyy]{datetime}
\usepackage{hhline}

\usepackage{listings}
\newcommand{\R}{\mathbb{R}}

\newcommand{\pr}{\mathbb{P}}
\newcommand{\F}{\mathscr{F}}

\DeclarePairedDelimiterX{\inner}[2]{\langle}{\rangle}{#1,#2}
\DeclarePairedDelimiterX{\norm}[1]{\|}{\|}{#1}

\DeclareMathOperator{\Tr}{Tr}

\theoremstyle{theorem}	\newtheorem{thm}{Theorem}[section]
\theoremstyle{theorem}	
\theoremstyle{theorem}	\newtheorem{prop}[thm]{Proposition}
\theoremstyle{theorem}	
\theoremstyle{definition}	 
\theoremstyle{definition}	\newtheorem{Def}[thm]{Definition}

\theoremstyle{theorem}	
\theoremstyle{definition}

\theoremstyle{plain}		
\theoremstyle{definition}	\newtheorem{rmk}[thm]{Remark}

\theoremstyle{definition}	\newtheorem{EX}{Example}[section]

\theoremstyle{definition}

\definecolor{OliveGreen}{cmyk}{0.64,0,0.95,0.40}
\definecolor{CadetBlue}{cmyk}{0.62,0.57,0.23,0}
\definecolor{lightlightgray}{gray}{0.93}

\begin{document}
\begin{center}
{\Large The Asymptotic Frequency of Stochastic Oscillators } \\
{\textsc{ (preprint)} }
\end{center}

\begin{center}
\begin{minipage}[t]{.75\textwidth}
\begin{center}
Zachary P.~Adams \\
Max Planck Institute for Mathematics in the Sciences\\
\hfill
\end{center}
\end{minipage}
\end{center}
\begin{center}
\today 
\end{center}
\vspace{.5cm}
\setlength{\unitlength}{1in}

\vskip.15in


{
  \centerline
  {\large \bfseries \scshape Abstract}
  \begin{quote}
We study stochastic perturbations of ODE with stable limit cycles -- referred to as \emph{stochastic oscillators} -- and investigate the response of the asymptotic (in time) frequency of oscillations to changing noise amplitude. 
Unlike previous studies, we do not restrict our attention to the small noise limit, and account for the fact that large deviation events may push the system out of its oscillatory regime. 
To do so, we consider stochastic oscillators conditioned on their remaining in an oscillatory regime for all time. 
This leads us to use the theory of \emph{quasi-ergodic measures}, and to define \emph{quasi-asymptotic frequencies}~as conditional, long-time average frequencies. 
We show that quasi-asymptotic frequencies exist under minimal assumptions, though they may or may not be observable in practice. 
Our discussion recovers and expands upon previous results on stochastic oscillators in the literature. 
In particular, existing results imply that the asymptotic frequency of a stochastic oscillator depends quadratically on the noise amplitude. 
We describe scenarios where this prediction holds, though we also show that it is not true in general -- potentially, even for small noise.  
\end{quote}
}


\section{Introduction}
\label{sec:Intro} 
\subsection{Background}
\label{subsec:Background} 
This paper is a small contribution to our understanding of how oscillatory dynamical systems respond to random perturbations. 
Specifically, we study ODE with stable periodic solutions, and investigate the effect which additive or multiplicative Guassian noise has on the asymptotic (time average) frequency of oscillations. 
When noise is present, we refer to these systems in general as \emph{stochastic oscillators}. 

We attempt to provide a unified, rigorous context for past results on the subject, and comment on several commonly made assumptions which do not hold in general. 
The mathematics employed in this paper is relatively simple. 
The only technical novelty of our approach is the use of ``quasi-ergodic measures'' (see Section \ref{sec:Unbounded}) to study the long term behaviour of stochastic oscillators. 

This subject -- the effect of noise on oscillatory dynamical systems -- has been studied from various perspectives over the past thirty years. 
In the 1990's, most work focused on the effect of noise on oscillatory dynamical systems' invariant measures. 
Particular attention was given to ``stochastic Hopf bifurcations'', for instance in the work of Arnold \&~Imkeller, \cite{AI98}, Arnold, Sri Namachchivaya, \&~Schenk-Hopp{\'e}~\cite{SH96}, and Baxendale \cite{B94}. 

A more dynamical approach gained popularity in the 2000's, based on earlier work of Guckenheimer \cite{G75}~and Winfree \cite{W74}~which rigorously defined the phase of deterministic oscillators. 
In this approach, a stochastic oscillator is projected, via a ``phase map'', onto the circle. 
The projected process is then identified as the ``phase'' of the system. 
There is some freedom in how one defines the phase, as we discuss below. 
Once the phase of an oscillatory dynamical system is defined, one may study how it, and properties derived from it, are affected by the addition of noise to the system. 

Following Teramae \&~Tanaka \cite{TT04}, the authors of \cite{B17}, \cite{BM18}, \cite{CQ21}, \cite{GPS18}, \cite{GTNE10}, \cite{TNE10}, \cite{YA07}, define phase maps by referring to the unperturbed (deterministic) oscillator. 
Meanwhile, \cite{C17}, \cite{CLT20}, \cite{EK21}, \cite{SP10a}, \cite{SP10b}, \cite{SPKR12}, \cite{SP13}, attempt to account for the non-deterministic behaviour of a stochastic oscillator when definining a phase map. 
Notably, the work of Schwabedal \&~Pikovsky in \cite{SP10a}, \cite{SP10b}, \cite{SPKR12}, \cite{SP13}, attempts to develop a theory of phase maps for noise-induced oscillations, rather than just oscillations perturbed by noise. 
We here note that the present paper focuses exlusively on noise-perturbed oscillations. 

As will be seen, the phase maps of Cao, Engel, Schwabedal, Pikovsky, and their collaborators in \cite{C17}, \cite{CLT20}, \cite{EK21}, \cite{SP10a}, \cite{SP10b}, \cite{SPKR12}, \cite{SP13}, do not strictly meet our definition of a phase map. 
However, the discussion of this paper could be generalized to include most of these definitions (except perhaps those of Schwabedal \&~Pikovsky in \cite{SP10a}, \cite{SP10b}) with little difficulty. 
See Remark \ref{rmk:Phase}~for further discussion of this point. 

In some of the works mentioned above, in particular \cite{B17}, \cite{BM18}, \cite{CQ21}, \cite{GPS18}, \cite{GTNE10}, \cite{TNE10}, \cite{TT04}, \cite{YA07}, approximate SDE for the phase of a stochastic oscillator are derived. 
From such an equation, \cite{B17}, \cite{GTNE10}, \cite{TNE10}, \cite{YA07}, derive a formula for a stochastic oscillator's asymptotic frequency. 
The results of each of these studies imply that  the difference between the asymptotic frequency of an oscillatory system with and without noise is proportional to the square of the noise amplitude. 
However, two common assumptions are made in their analyses, which do not hold in general. 

First, the authors assume that stochastic deviations in the amplitude of oscillations due to noisy perturbations are negligible. 
It is true that large deviation results \cite{FW98}~guarantee that the system should stay in a ``small'' neighbourhood of a stable deterministic limit cycle for a ``long time'' (in a sense which we do not make precise here). 
However, if a deterministic limit cycle has a bounded basin of attraction, and the system is perturbed by additive noise (as used in several of the references we cite), then at some almost surely finite stopping time the stochastic system will exit the basin of attraction of the limit cycle. 
From this time onwards, the phase of the system is not well-defined, and hence neither is the time average of the system's frequency. 

Second, it is assumed that the pointwise ergodic theorem\footnote{
For details of the pointwise ergodic theorem which we use in Section \ref{sec:Bounded}~(first proven by Birkhoff \cite{B31}), refer to Section 1.2 of the textbook by Cornfeld \emph{et al.}~\cite{CFS12}. 
}
can be applied to the phase of a stochastic oscillator. 
In most cases this is true -- and would always be true if the assumption of the previous paragraph held -- but it is a technicality which should be taken care of. 

Not all studies of stochastic oscillators make these assumptions. 
In Bressloff \&~MacLaurin \cite{BM18}, Cheng \&~Qian \cite{CQ21}, and Giacomin \emph{et al.}~\cite{GPS18}, SDE for the phase of a stochastic oscillator are rigorously derived.  
The authors are careful to note that the phase may only be well-defined up to some finite stopping time. 

While the analysis of \cite{BM18}~ends at a finite stopping time, \cite{CQ21}, \cite{GPS18}, consider long time dynamics by making use of large deviation theory. 
In \cite{CQ21}, \cite{GPS18}, time and noise amplitude are simultaneously taken to infinity and zero, respectively, in a co-dependent fashion. 
The phase is thereby guaranteed to be well-defined throughout their arguments. 
Similar to \cite{B17}, \cite{GTNE10}, \cite{TNE10}, \cite{YA07}, the results of \cite{GPS18}~imply that the difference between the average frequency of an oscillatory system with and without noise is proportional to the square of the noise amplitude in the small noise regime.  

In the present study, we first consider (in Section \ref{sec:Bounded}) the long term dynamics of stochastic oscillators by assuming that the phase is always well-defined. 
Then, dropping this assumption in Section \ref{sec:Unbounded}, we use the theory of quasi-ergodic measures (introductions to which can be found in the textbook of Collet \emph{et al.}~\cite{CMSM12}~or the recent dissertation of Villemonais \cite{V19}), which allows us to consider a stochastic oscillator conditioned on the phase being well-defined for all time. 
We thence come to a formula for a ``quasi-asymptotic frequency''. 
We see a clear possibility of non-quadratic dependence of this quasi-asymptotic frequency on the noise amplitude. 
Our formula for the (quasi) asymptotic frequency of a stochastic oscillator is compared with previous results in Section \ref{sec:Comparison}. 
In particular, our simple arguments recover the results of Giacomin \emph{et al.}~\cite{GPS18}~up to a correction term. 
We conjecture that this correction term may be non-trivial in the small-noise regime for certain systems, though further investigation is needed to prove this. 
When our prediction differs from that of \cite{GPS18}, the asymptotic frequencies which we predict are likely observable on different time scales from those identified in \cite{GPS18}. 

In Section \ref{sec:Conclusions}~we discuss the advantages and disadvantages of our work. 
As we will see, the primary advantages are that it is rigorous, and that it is not restricted to the small-noise amplitude regime. 
A major disadvantage is that, though quasi-asymptotic frequencies are always well-defined, they may not be ``observable'', in a sense made precise below. 
Finally, directions for future work are considered.

\subsection{Technical setup}
\label{sec:Setup}
Before proceeding, we introduce our technical setup and state necessary definitions. 
Past results and present goals are also stated in a more precise language than in Section \ref{subsec:Background}. 

We begin by considering an autonomous differential equation generated by a $C^2$ vector field $V$ on $\R^d$, 
\begin{equation}
\label{eq:ODE}
\partial_tx=V(x). 
\end{equation} 
Denote the flow map of \eqref{eq:ODE}~by $(t,x)\mapsto X_t^0(x)$. 
We assume throughout this document that \eqref{eq:ODE}~has a stable limit cycle $\Gamma$ of period $T>0$, in the sense of \cite{G75}. 
That is, $\Gamma$ is a one dimensional manifold parameterized as $\Gamma=\{\gamma_t\}_{t\in\R}$ such that 
\begin{enumerate}
\item $X_t^0(\gamma_s)\,=\,\gamma_{t+s}$ and $\gamma_{t+T}=\gamma_t$, and 
\item there exists an open set $B(\Gamma)\subset\R^d$ containing $\Gamma$ such that for all $x\in B(\Gamma)$, 
\[
\lim_{t\rightarrow\infty}\inf_{y\in\Gamma}\norm*{X_t^0(x)-y}\,\xrightarrow[t\rightarrow\infty]{}\,0. 
\]
\end{enumerate} 

Let $B(\Gamma)$ be the basin of attraction of $\Gamma$. 
We take a surjective map $\pi:B(\Gamma)\rightarrow[0,T)$ such that 
\begin{equation}
\label{eq:Phase}
\pi(\gamma_s)\,=\,s\,\,\text{mod}\,T, 
\end{equation}
and let the \emph{$\pi$-phase}~(or simply \emph{phase}, when unambiguous) of a continuous path $(Y_t)_{t\ge0}$ in $B(\Gamma)$ be $\pi(Y_t)$.  
In general, we refer to such a map as a \emph{phase map}. 
For the rest of the document, all phase maps are assumed to be $C^2$. 
For a phase map $\pi$, the corresponding \emph{normalized phase map}~is $\pi_1\coloneqq T^{-1}\pi$. 

\begin{rmk}
\label{rmk:Phase} 
The ``phase maps'' of Cao \cite{C17}, Cao \emph{et al.}~\cite{CLT20}, Engel \&~Kuehn \cite{EK21}, and Schwabedal \&~collaborators \cite{SP10a}, \cite{SP10b}, \cite{SPKR12}, \cite{SP13}, do not satisfy condition \eqref{eq:Phase}. 
Dropping this condition, the discussion of this paper would extend to the phase maps of  \cite{C17}, \cite{CLT20}, \cite{SP13}. 
For our discussion to apply to the phase map of \cite{EK21}, we would need to consider random phase maps. 
The complications introduced to our arguments by weakening the definition of a phase map are marginal. 
Nevertheless, we maintain the condition \eqref{eq:Phase}~in the definition of a phase map for didactic purposes. 
Our arguments do not easily extend to the phase maps of \cite{SP10a}, \cite{SP10b}, which are designed to handle noise-induced (rather than noise-perturbed) oscillations. Nor do our arguments easily extend to the phase map used in \cite{SPKR12}, which deals with chaotic (rather than stochastic) oscillators. 
\end{rmk}

In the literature, frequent use is made of the \emph{isochron map}, which is the unique phase map invariant under the time $T$ flow of \eqref{eq:ODE}. 
That is, the isochron map is the unique map $\pi:B(\Gamma)\rightarrow[0,T)$ satisfying \eqref{eq:Phase}~and 
\begin{equation}
\label{eq:MonodromeInvariance}
\pi(x)=\pi(X_T^0(x))\,\,\,\,\text{ for all }\,x\in B(\Gamma). 
\end{equation}
Equivalently, the isochron map can be defined for each $x\in B(\Gamma)$ as the unique number $\pi(x)\in[0,T)$ such that 
\[
\norm*{X_t^0(x) - \gamma_{t+\pi(x)}}\,\xrightarrow[t\rightarrow\infty]{}\,0. 
\]
The isochron map was introduced in the context of deterministic oscillators by Winfree \cite{W74}, and studied from a mathematical perspective by Guckenheimer \cite{G75}~(though as \cite{G75}~points out, similar ideas can be traced back to Poincar{\'e}). 
When $\pi$ is the isochron map, we refer to the $\pi$-phase of a continuous path in $B(\Gamma)$ as its \emph{isochronal phase}. 

Throughout our discussion, we must take care of possible singularities of phase maps $\pi$ that may exist at the boundary of $B(\Gamma)$. 
A \emph{phase singularity}~is a point $x_0\in\partial B(\Gamma)$ such that some derivative of $\pi$ blows up at $x_0$. 
For instance, if $\pi$ is the isochron map and $x_0$ is a zero of $V$ with an unstable manifold intersecting $B(\Gamma)$, then $x_0$ is a phase singularity of $\pi$. 
This follows from the fact that $\pi'(x)V(x)=1$ for all $x\in B(\Gamma)$, so 
\[
\norm*{\pi'(x)}\,\sim\,\norm*{V(x)}^{-1}\,\xrightarrow[x\rightarrow x_0]{}\,\infty. 
\]

Our interest lies in stochastic perturbations of \eqref{eq:ODE}. 
Let $(\Omega,\F,(\F_t)_{t\ge0},\pr)$ be a filtered probability space satisfying the usual conditions on which all random variables and stochastic processes of this document are to be defined. 
Consider the stochastic differential equation 
\begin{equation}
\label{eq:SDE}
dX=V(X)\,dt + \sigma B(X)\,dW, 
\end{equation}
where $W=(W_t)_{t\ge0}$ is a standard $(\F_t)_{t\ge0}$-adapted Brownian motion on $\R^d$, $B:\R^d\rightarrow\R^{d\times d}$ is a continuous map, and $\sigma\ge0$. 
A \emph{solution}~to \eqref{eq:SDE}~is a continuous, $(\mathcal{F}_t)_{t\ge0}$ adapted stochastic process $(X_t)_{t\ge0}$ 
such that 
\begin{equation}
\label{eq:IntSDE}
X_t\,=\,X_0 + \int_0^tV(X_s)\,ds + \sigma\int_0^tB(X_s)\,dW_s\quad\text{ almost surely}. 
\end{equation}
For details on the solution theory of \eqref{eq:SDE}, refer to \cite{O10}. 

For $x\in\R^d$, let $\pr_x(\,\cdot\,)\coloneqq\pr\left(\,\cdot\,|X_0=x\right)$, where $X_0$ is the initial condition of \eqref{eq:SDE}. 
For any distribution $\nu$ on $\R^d$ let 
\[
\pr_\nu\left(\,\cdot\right)\,\coloneqq\,\int_{\R^d}\pr_x(\,\cdot\,)\,\nu(dx), 
\]
which is $\pr$ conditioned on the initial distribution of \eqref{eq:SDE}~being $\nu$. 
Denote by $\mathbb{E}_x$ and $\mathbb{E}_\nu$ the expectation with respect to $\pr_x$ and $\pr_\nu$. 

Let the flow map of \eqref{eq:SDE}~be $(t,x)\mapsto X_t^\sigma(x)$. 
When unambiguous, we omit the dependence of the flow on initial conditions, writing $X_t^\sigma=X_t^\sigma(x)$. 
For any phase map $\pi$, note that the phase of \eqref{eq:SDE}~is only defined so long as $X_t^\sigma$ remains in $B(\Gamma)$ for all $t\ge0$, which is not guaranteed in general. 
The \emph{exit time}~of $(X_t^\sigma)_{t\ge0}$ from $B(\Gamma)$ is 
\[
\tau_\sigma\,\coloneqq\, \inf\left\{t>0\,:\,X_t^\sigma\in\partial B(\Gamma)\right\}, 
\]
so that any phase map of \eqref{eq:SDE}~is only defined for $t<\tau_\sigma$. 

The main object in which we are interested is the \emph{asymptotic frequency}~of \eqref{eq:SDE}, defined as 
\[
c_\sigma\,\coloneqq\,\lim_{t\rightarrow\infty}\frac{1}{t}\pi_1(X_t^\sigma) 
\]
when it exists. 
We use the normalized phase map, so that $c_\sigma$ is the long-time average number of full rotations per unit time. 
If the unnormalized phase map were used (as in \cite{GPS18}), then $c_\sigma$ would be the long-time average number of full rotations per deterministic period. 

In Section \ref{sec:Bounded}~we assume that $\tau_\sigma=\infty$, in which case $c_\sigma\in\R$ is (usually) well-defined. 
A formula for the asymptotic frequency is given, which facilitates a qualitative understanding of its dependence on the noise amplitude $\sigma>0$. 
In Section \ref{sec:Unbounded}~we allow for $\tau_\sigma<\infty$ almost surely, and study $c_\sigma$ conditioned on the event $\tau_\sigma=\infty$ using the theory of quasi-ergodic measures. 
We discuss examples where $c_\sigma$ is observable, in the sense that $t^{-1}\pi_1(X_t^\sigma)$ approaches $c_\sigma$ for some $t<\tau_\sigma$ with high probability. 
In Section \ref{sec:Comparison}, we study the qualitative dependence of $c_\sigma-c_0$ on $\sigma>0$, and compare our results with those of past studies, 
\cite{GPS18}~in particular.  
Section \ref{sec:Conclusions}~concludes the paper, summarizing our results and highlighting directions for future research.

\section{Systems which oscillate for all time } 
\label{sec:Bounded} 
In this section, we study the asymptotic frequency of a stochastic oscillator assuming that $\tau_\sigma=\infty$. 
Note that in our setup, this is usually only possible for multiplicative noise, specifically when the diffusion coefficient in\eqref{eq:SDE}~is such that $B(x)\rightarrow0$ as $x\rightarrow\partial B$. 

Our main result is the following theorem; its analogue when $\tau_\sigma$ is almost surely finite is given in Theorem \ref{thm:QuasiFrequency}. 
As noted in Remark \ref{rmk:Unbounded}, the assumption that $\mu_\sigma$ has bounded support in $B(\Gamma)$ can be substantially weakened.
We include this assumption at first to simplify the theorem's proof. 

\begin{thm}
\label{thm:ErgodicFrequency}
Let $(X_t^\sigma)_{t\ge0}$ be the stochastic process governed by \eqref{eq:SDE}, and assume $\tau_\sigma=\infty$. 
Let $\Gamma$ be a stable limit cycle of \eqref{eq:ODE}~with basin of attraction $B(\Gamma)$ and period $T>0$. 
Fix a phase map $\pi:B(\Gamma)\rightarrow[0,T)$. 
Suppose that 
\begin{enumerate}[(i)]
\item $(X_t^\sigma)_{t\ge0}$ has an ergodic measure $\mu_\sigma$ with bounded support in $B(\Gamma)$, 
\item If $\pi$ has a phase singularity $x_0\in\partial B(\Gamma)$, then $\mu_\sigma$ is such that for $\delta>0$ 
\begin{equation}
\label{eq:Integrable}
\int_{B_\delta(x_0)} \norm*{\pi'(x)}\,\mu_\sigma(dx)\,<\,\infty,\quad \int_{B_\delta(x_0)} \norm*{\pi''(x)}\,\mu_\sigma(dx)\,<\,\infty, 
\end{equation}
\emph{i.e}~$\mu_\sigma$ decays sufficiently fast near $x_0$. 
\end{enumerate}
Then, 
\begin{equation}
\label{eq:ErgodicFrequency}
\begin{aligned}
c_\sigma\,&\coloneqq\,\lim_{t\rightarrow\infty}\frac{1}{t}\pi_1(X_t^\sigma)\\
&=\,\int_{B(\Gamma)}\pi_1'(x)V(x) + \frac{\sigma^2}{2}\Tr\pi_1''(x)[B(x),B(x)]\,\mu_\sigma(dx) 
\end{aligned}
\end{equation}
exists as a deterministic real number, the limit converging in probability. 
\begin{proof}
By It{\^o}'s formula, we have 
\begin{equation}
\label{eq:piIto}
\begin{aligned}
\pi(X^\sigma_t)\,&=\, \int_0^t \pi'(X^\sigma_s)V(X^\sigma_s) + \frac{\sigma^2}{2}\Tr\pi''(X^\sigma_s)[B(X^\sigma_s)\,\cdot\,,B(X^\sigma_s)\,\cdot\,]\,ds \\
&\qquad+ \sigma\int_0^t\pi'(X^\sigma_s)B(X^\sigma_s)\,dW_s \\
&\eqqcolon\, I_t + II_t . 
\end{aligned}
\end{equation}
To handle $I_t$, we remark that condition (ii), the assumed regularity of $\pi$, and the boundedness of the support of $\mu_\sigma$ imply that the map 
\begin{equation}
\label{eq:Lpi}
x\,\mapsto\,\pi'(x)V(x) + \frac{\sigma^2}{2}\Tr\pi''(x)[B(x)\,\cdot\,,B(x)\,\cdot\,] 
\end{equation}
is integrable over $B(\Gamma)$ with respect to $\mu_\sigma$. 
Applying the pointwise ergodic theorem (see Section 1.2 of \cite{CFS12}) with respect to $\mu_\sigma$ then yields 
\[
\begin{aligned}
&\lim_{t\rightarrow\infty}\frac{1}{t}\int_0^t\pi'(X^\sigma_s)V(X^\sigma_s) + \frac{\sigma^2}{2}\Tr\pi''(X^\sigma_s)[B(X^\sigma_s)\,\cdot\,,B(X^\sigma_s)\,\cdot\,]\,ds \\
&\qquad=\, \int_{\Gamma_\delta} \pi'(x)V(x) + \frac{\sigma^2}{2}\Tr\pi''(x)[B(x)\,\cdot\,,B(x)\,\cdot\,]\,\mu_\sigma(dx).  
\end{aligned}
\]

For $II_t$, assume that the initial distribution of $(X_t)_{t\ge0}$ is $\mu_\sigma$. 
By It{\^o}'s isometry and the definition of the distribution of a process, 
\[
\begin{aligned}
\mathbb{E}_{\mu_\sigma}\left[\norm*{II_t}^2\right]\,&\le\, \mathbb{E}_{\mu_\sigma}\left[\int_0^t\norm*{\pi'(X_s)B(X_s)}^2\,ds\right]\\
&\le\, \int_{B(\Gamma)}\int_0^t\norm*{\pi'(x)B(x)}^2\,ds\,\mu_\sigma(dx) \\
&\le\, C_0 t. 
\end{aligned}
\]
By the Burkholder-Davis-Gundy inequality, 
\[
\mathbb{E}_{\mu_\sigma}\left[[II]_t\right]\,\le\,C_1t 
\]
for some $C_1>0$, where $[II]_t$ denotes the quadratic variation of $II_t$. 
This implies that $\mathbb{E}_{\mu_\sigma}\left[t^{-2}[II]_t\right]\rightarrow0$ as $t\rightarrow\infty$. 
Then, by Theorem 4.1 of van Zanten \cite{vZ00}, 
\[
\frac{1}{t}II_t\,\xrightarrow[t\rightarrow\infty]{}\,0\qquad\text{ in distribution.} 
\]
As the distributional limit of $t^{-1}II_t$ is zero, a constant, this is also a limit in probability. 
If the initial distribution of $(X_t)_{t\ge0}$ is not $\mu_\sigma$, we nevertheless know that the distribution of $X_t$ becomes arbitrarily close to $\mu_\sigma$ as $t\rightarrow\infty$. 
Therefore, we apply the same argument after waiting for transient behaviour to die out. 
\end{proof}
\end{thm}

When $c_\sigma$ exists, as defined in \eqref{eq:ErgodicFrequency}, we refer to it as the \emph{asymptotic frequency}~of \eqref{eq:SDE}~in $B(\Gamma)$.  

Theorem \ref{thm:ErgodicFrequency}~is really just an application of the pointwise ergodic theorem. 
We are not the first to use the pointwise ergodic theorem to study the asymptotic frequency of a stochastic oscillator, see for instance Teramae \emph{et al.}~\cite{TNE10}~and Yoshimura \&~Arai \cite{YA07}. 
However, \cite{YA07}~assumes that the distance of $X_t^\sigma$ from $\Gamma$ is approximately zero for all $t>0$, while \cite{TNE10}~only considers the small noise limit of the system and assumes that attraction of \eqref{eq:ODE}~to $\Gamma$ occurs infinitely fast. 
Therefore, the integrability condition (ii) does not enter into their consideration. 
As already remarked, the assumption that $X_t^\sigma$ remains near $\Gamma$ for all $t>0$ is not necessarily good, even for small $\sigma>0$. 

When attempting to apply Theorem \ref{thm:ErgodicFrequency}, the integrability condition \eqref{eq:Integrable}~needs to be checked. In Appendix \ref{app:Sufficient}, we provide sufficient conditions for \eqref{eq:Integrable}~to be satisfied. 
In Example \ref{EX:LambdaOmega}, we consider a simple SDE satisfying this condition for some values of $\sigma>0$. 
However, Examples \ref{EX:LambdaOmega}~\&~\ref{EX:LambdaOmegaAdditive}~illustrate that Appendix \ref{app:Sufficient}~provides sufficient, but not necessary, conditions for the existence of an asymptotic frequency. 

\begin{rmk}
\label{rmk:ConvergenceRates}
We have not been able to obtain estimates on the rate of convergence in \eqref{eq:ErgodicFrequency}, though in some cases we expect this rate to be exponential. 
Indeed, if $(X_t)_{t\ge0}$ is restricted to a compact subset of its phase space, it satisfies Doeblin's condition (see for instance Section 16 of Meyn \&~Tweedie \cite{MT12}). 
Moreover, if $\pi$ does not have a phase singularity on $\partial B(\Gamma)$, then the integrand appearing in $I_t$ is a bounded function on $B(\Gamma)$. 
Hence, for $I_t$ the conditions of Katz \&~Thomasian \cite{KT61}~are satisfied, and we have that $\frac{1}{t}I_t$ converges to $c_\sigma$ at an exponential rate. 

However, this exponential rate of convergence is not guaranteed in the presence of a phase singularity on $\partial B(\Gamma)$. 
Nor is an exponential rate of convergence of $\frac{1}{t}II_t$ to zero guaranteed by any theory which we are aware of. 
Therefore, we cannot conclude that the rate of convergence of $\frac{1}{t}\pi_1(X_t)$ to $c_\sigma$ is exponential. 
Future work may study this rate of convergence, and rates of convergence to quasi-ergodic averages in general. 
\end{rmk}

\begin{rmk}
\label{rmk:Unbounded}
In condition (i) of Theorem \ref{thm:ErgodicFrequency}, we require the support of $\mu_\sigma$ to be bounded. 
This is to prevent the potential growth of the drift and diffusion coefficients $V,\,B,$ at infinity from preventing the $\mu_\sigma$-integrability of the functional in \eqref{eq:Lpi}.  
However, the condition is not necessary. 
Indeed, suppose that 
\begin{itemize} 
\item $B(\Gamma)=\R^d$, so that the phase is well-defined for all time, and that 
\item the distribution of $(X_t)_{t\ge0}$ converges to $\mu_\sigma$ at an exponential rate, which is the case for a wide range of processes. 
\end{itemize}
Then, $\mu_\sigma$ satisfies a Poincar{\'e}~inequality, which implies that $\mu_\sigma$ must have tails which decay to zero at an exponential rate (see Chapter 4 of Bakry \emph{et al.}~\cite{BGL14}~for details). 
So long as $V$ and $B$ have subexponential growth as $\norm*{x}\rightarrow\infty$, we may therefore conclude that the functional in \eqref{eq:Lpi}~is still integrable with respect to $\mu_\sigma$, and the argument of Theorem \ref{thm:ErgodicFrequency}~applies to this scenario with little extra effort. 
\end{rmk}

\begin{figure}[h!]
\centering
\begin{subfigure}[b]{0.45\textwidth}
\centering
\includegraphics[width=\textwidth]{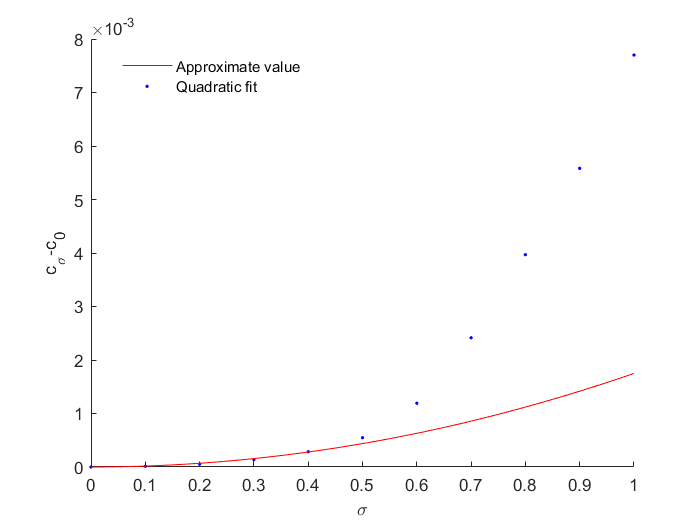}
\caption{}
\label{fig:LambdaOmega}
\end{subfigure}
\hfill
\begin{subfigure}[b]{0.45\textwidth}
\centering
\includegraphics[width=\textwidth]{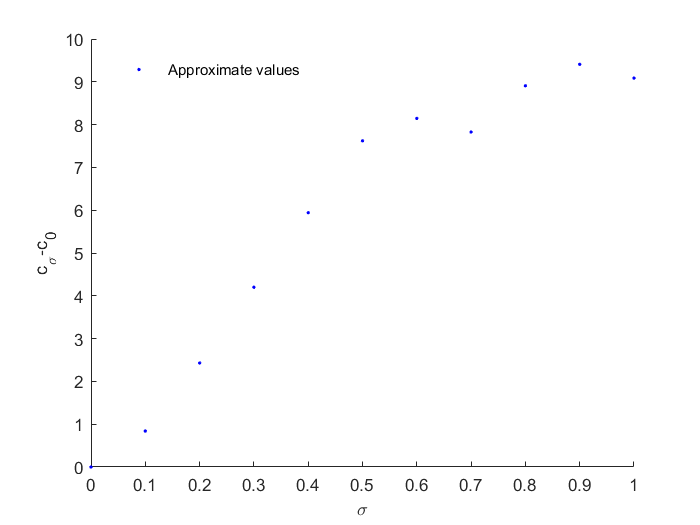}
\caption{}
\label{fig:LambdaOmegaAdditive}
\end{subfigure}
\caption{
Approximate values of $c_\sigma-c_0$ for $\sigma\in[0,10]$ in \eqref{eq:LambdaOmegaBounded}~\&~\eqref{eq:LambdaOmegaBoundedUnsym}, shown in (a) and (b) respectively. 
Both plots were obtained by simulating the respective system using an Euler-Maruyama scheme up to time $t=1000$, at which point convergence appears to be reached. Note that the quadratic fit \eqref{eq:cApprox}~appears to be good in (a) for $\sigma\in[0,10]$, while super-quadratic dependence appears to occur for $\sigma\gtrsim5$ in (b).   
}
\label{fig:Frequencies}
\end{figure}

\begin{EX}
\label{EX:LambdaOmega}
Consider the following stochastic perturbation of the Hopf normal form, 
\begin{equation}
\label{eq:LambdaOmegaBounded}
\begin{aligned}
dx\,&=\,  (x - y -x(x^2+y^2) )\,dt + \sigma x(2-(x^2+y^2))\,dW,\\
dy\,&=\, (x + y - y(x^2+y^2) )\,dt + \sigma y(2-(x^2+y^2))\,dW, 
\end{aligned}
\end{equation}
driven by a single Brownian motion $(W_t)_{t\ge0}$. 
Here and in all other examples, the noise is interpreted in the It{\^o}~sense. 
Note that the choice of noise in \eqref{eq:LambdaOmegaBounded}~implies that its solution is bounded in the circle of radius $\sqrt{2}$ uniformly in time. 

When $\sigma=0$, the unit circle $\Gamma=S^1$ is globally stable (\emph{i.e.}~$B(\Gamma)=\R^2$).  
Hence $\tau_\sigma=\infty$, regardless of the choice of diffusion coefficient, so that the isochronal phase of the solution process is defined for all $t\ge0$. 
This example is particularly nice for our purposes, since one may check that the isochron of each $(x,y)\in S^1$ is the ray 
\[
\pi^{-1}(x,y)\,=\,\big\{s(x,y)\,:\,s>0\big\}. 
\]

Note that the diffusion coefficient in \eqref{eq:LambdaOmegaBounded}~is such that the conditions of Corollary \ref{cor:}~are met, so long as 
\[
0\,\le\,\sigma\,\le\, \sigma_*\,\coloneqq\,\Tr V'(0)/2d\,=\,1/2. 
\] 
Thus, $t^{-1}\pi(X_t)$ converges to some $c_\sigma\in\R$ in probability as $t\rightarrow\infty$ for $\sigma<\sigma_*$. 
However, our numerical experimets suggest that this convergence also occurs for some $\sigma>\sigma_*$. 

In Figure \ref{fig:LambdaOmega}, we approximate $c_\sigma$ by a numerical value of $t^{-1}\pi(X_t^\sigma)$ for large $t>0$. 
We ran an Euler-Maruyama scheme (see Chapter 8 of \cite{LPS14}) with a time step of $dt = 0.0025$ up to time $t_{end}=5000$. 
Note that the approximation 
\begin{equation}
\label{eq:cApprox}
c_\sigma-c_0\,\simeq\, m\sigma^2 \quad\text{ with }\quad m\,\coloneqq\,1.8*10^{-3}. 
\end{equation}
appears to be good for $\sigma\in[0,0.4]$. 
Thus in this example, the prediction of Giacomin \emph{et al.}~\cite{GPS18}~appears to extend into the moderate-noise regime, and holds on time-scales much larger than what they predict. 
For $\sigma>0.4$, $c_\sigma$ appears to enter a non-quadratic regime. 

As noted in Remark \ref{rmk:Unbounded}, uniform in time boundedness of the solution is a sufficient condition for the existence of a deterministic asymptotic frequency, but not necessary. 
For instance, consider 
\begin{equation}
\label{eq:LambdaOmega}
\begin{aligned}
dx\,&=\,  (x - y -x(x^2+y^2) )\,dt + \sigma x\,dW,\\
dy\,&=\, (x + y - y(x^2+y^2) )\,dt + \sigma y\,dW. 
\end{aligned}
\end{equation}
With noise of this form the process is guaranteed to leave any bounded subset of $\R^2$ at some almost surely finite stopping time. 
However, the system still possesses a unique ergodic measure in $B(\Gamma) = \R^2/\{0\}$, and satisfies condition (ii) of Theorem \ref{thm:ErgodicFrequency}~for $\sigma<\sigma_*$, by Corollary \ref{cor:}. 
Numerical simulations again suggest that the asymptotic frequency exists, and depends quadratically on $\sigma\in[0,0.5]$ (not shown). 
\end{EX}

\begin{figure}[h!]
\centering
\begin{subfigure}[b]{0.45\textwidth}
\centering
\includegraphics[width=\textwidth]{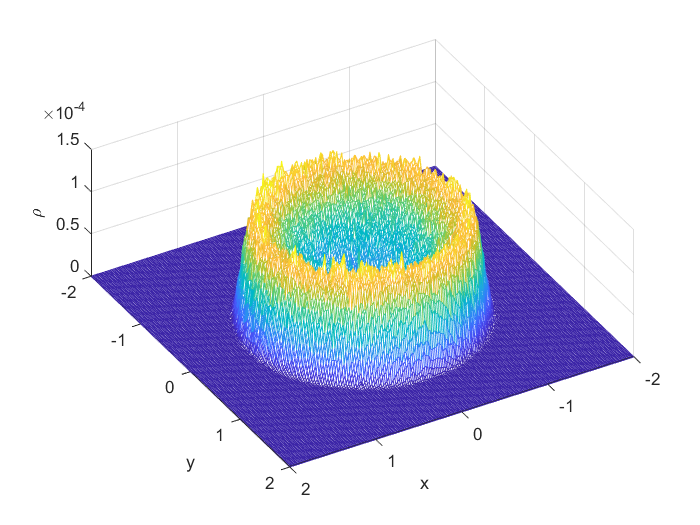}
\caption{}
\label{subfig:LOInvariant}
\end{subfigure}
\begin{subfigure}[b]{0.45\textwidth}
\centering
\includegraphics[width=\textwidth]{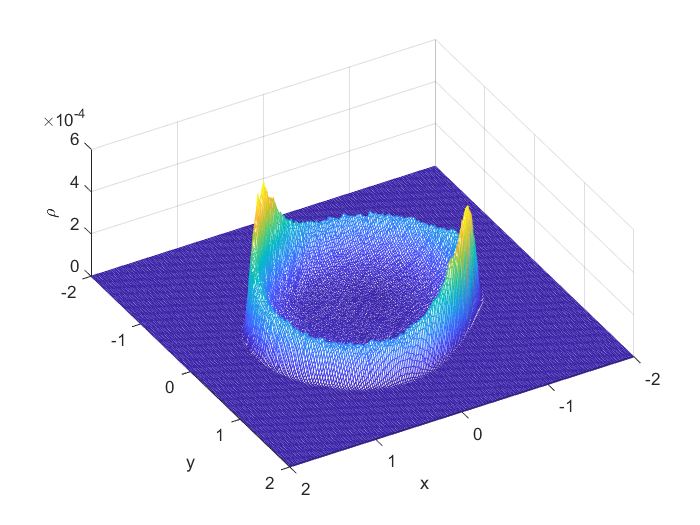}
\caption{}
\label{subfig:LOAInvariant}
\end{subfigure}
\caption{ 
Monte-Carlo approximations of the unique ergodic measures in $B(\Gamma)$ of \eqref{eq:LambdaOmegaBounded}~and \eqref{eq:LambdaOmegaBoundedUnsym}~(shown in (a) and (b), respectively), both with $\sigma=0.4$. 
In (a), note the approximate radial symmetry of the invariant measure, not present in (b). 
 }
\label{fig:LambdaOmegaInvariantMeasures}
\end{figure}

\begin{EX}
\label{EX:LambdaOmegaAdditive} 
We consider a process generated by an SDE to which Corollary \ref{cor:}~does not apply, 
\begin{equation}
\label{eq:LambdaOmegaBoundedUnsym}
\begin{aligned}
dx\,&=\,  (x - y -x(x^2+y^2) )\,dt + \sigma (2-(x^2+y^2))\,dW,\\
dy\,&=\, (x + y - y(x^2+y^2) )\,dt + \sigma (2-(x^2+y^2))\,dW. 
\end{aligned}
\end{equation}
Investigating the Fokker-Planck operator of this system, one finds that the comparison principle cannot be used to conclude that its ergodic measure (which exists and is unique) decays quadratically near zero. 
Hence, our previous arguments do not guarantee that the integral in \eqref{eq:ErgodicFrequency}~converges. 

Nevertheless, simulating the system for a long time, we observe the apparent convergence of $t^{-1}\pi(X_t)$ to a deterministic $c_\sigma$ for $\sigma\in[0,0.5]$. 
The asymptotic frequency appears to depend non-quadratically on $\sigma\in[0,0.5]$. 
Beyond the value $\sigma=0.5$, we did not observe apparent convergence of $\frac{1}{t}\pi_1(X_t)$ to a fixed deterministic constant, though this may only be due to an exceptionally slow rate of convergence. 
See Figure \ref{subfig:LOAInvariant}. 

As we discuss further in Section \ref{sec:Comparison}, Theorem \ref{thm:ErgodicFrequency}~implies that a non-quadratic dependence of $c_\sigma$ on $\sigma>0$ must be due to significant changes in $\mu_\sigma$. 
Indeed, in this example, the radially asymmetric noise which we use leads to a loss of radial symmetry in $\mu_\sigma$. 
See Figure \ref{subfig:LOAInvariant}. 
We may further remark that the asymmetry of $\mu_\sigma$ in this example corresponds to a sort of noise-induced bistability in \eqref{eq:LambdaOmegaBoundedUnsym}. 
Each of the peaks of the invariant measure seen in Figure \ref{subfig:LOAInvariant}~corresponds to a meta-stable state, between which the system rapidly switches. 
This sort of noise-induced bistability has also been studied by Newby \&~Schwemmer \cite{NS14}, \cite{SN15}. 
\end{EX}

\begin{EX}
\label{EX:FHN} 
The theory developed here extends, in some cases, to parabolic SPDE interpreted as stochastic evolution equations on a Hilbert space. 
For instance, consider 
\begin{equation}
\label{eq:SPDE}
\begin{aligned}
dX\,=\,\big(\Delta X + N(X)\big)\,dt + \sigma B(X)\,dW, 
\end{aligned}
\end{equation}
with periodic spatial domain $S^1$, interpreted as an evolution equation on $L^2(S^1;\R^n)$. 
The Laplace operator is $\Delta:H^1(S^1;\R^d)\subset L^2(S^1;\R^d)\rightarrow L^2(S^1;\R^d)$, and $N:C_b(S^1;\R^d)\rightarrow C_b(S^1;\R^d)$ is a locally Lipschitz nonlinearity. 
For simplicity, we assume that $B(x)$ is a trace class operator for each $x\in L^2(S^1;\R^d)$, that $x\mapsto B(x)$ is continuous in the strong operator topology, and that $(W_t)_{t\ge0}$ is a trace-class Wiener process. 
These assumptions allow us to conclude that existence and uniqueness of solutions to \eqref{eq:SPDE}~hold in the mild sense. 
That is, there exists a unique stochastic process $(X_t^\sigma)_{t\ge0}$ satisfying 
\[
X_t^\sigma\,=\,e^{t\Delta}X_0^\sigma + \int_0^t e^{\Delta(t-s)}N(X_s^\sigma)\,ds + \sigma\int_0^t e^{\Delta(t-s)} B(X_s^\sigma)\,dW_s 
\]
for $X_0^\sigma\in C_b(S^1;\R^d)$. 
See Da Prato \&~Zabczyk \cite{DPZ14}~or Liu \&~R{\"o}ckner \cite{LR16}~for details. 

When $\sigma=0$, there are instances of \eqref{eq:SPDE}~possessing an asymptotically stable travelling wave solution in $L^2(S^1;\R^d)$, for instance the FitzHugh-Nagumo system  
\begin{equation}
\label{eq:FHN}
\begin{aligned}
du\,&=\, \big(D\Delta u + u(1-u)(u-a)\big)\,dt + \sigma  u(1-u)(u-a)\,dW, \\
dv\,&=\, \big(\delta D\Delta v + \gamma u - v\big)\,dt. 
\end{aligned}
\end{equation}
When $\sigma=0$, the existence and stability of a travelling wave solution of \eqref{eq:FHN}~on a periodic spatial domain is proven in Ariola \&~Koch \cite{AK16}, as long as the scaling parameter $D$ is sufficiently small (equivalently, as long as the spatial domain is sufficiently large). 
The choice of diffusion coefficient in \eqref{eq:FHN}~guarantees that the solution remains in a bounded subset of the travelling wave's basin of attraction for all time. 
On a periodic spatial domain, the travelling wave is a periodic solution. 
When the noise is trace class, MacLaurin \cite{M21}~guarantees that the isochron map is well-defined, $C^2$, and satisfies an It{\^o}~formula. 
This allows us to translate most of the proof of Theorem \ref{thm:ErgodicFrequency}~to this setting unchanged. 

So long as $\delta>0$, we may apply Theorem 11.38 of Da Prato \&~Zabczyk \cite{DPZ14}~to conclude that \eqref{eq:FHN}~possesses an ergodic measure in $B(\Gamma)$. 
Theorem \ref{thm:ErgodicFrequency}~therefore translates to this setting, and we can conclude that \eqref{eq:FHN}~has an asymptotic frequency in $B(\Gamma)$ for all $\sigma>0$. 
This asymptotic frequency is in fact the asymptotic speed of the stochastic travelling wave. 

Future work will focus on studying this system with additive noise, using the theory outlined in the following section. 
We are also interested in extending this example to the case of non-trace class noise, and the case $\delta=0$. 
When the noise is not trace class, the results of \cite{M21}~need to be strengthened. 
When $\delta=0$, we must handle the fact that the linear part of the system's drift coefficient does not generate a compact $C_0$-semigroup. 
This makes proving the existence of an ergodic measure slightly more difficult. 

One might also consider how the speed of the wave on an unbounded spatial domain compares to the speed of the wave on a ``large'' periodic spatial domain. 
This would provide a rigorous foundation for the existence of the asymptotic stochastic wave speeds computed \emph{e.g.}~in the thesis work of Hamster \cite{H19}, \cite{HH20b}. 
It could also provide a theoretical framework for other works on the effects of noise on the speed of travelling waves, such as Eichinger \emph{et al.}~\cite{EGK20}~or MacLaurin \cite{M21}. 
\end{EX}

\section{Systems which may not oscillate for all time}
\label{sec:Unbounded} 
In many situations, one cannot assume that $\tau_\sigma=\infty$. 
Nevertheless, we often expect the solution of \eqref{eq:SDE}~to remain in $B(\Gamma)$ for a ``long time''. 
For instance, persistence in $B(\Gamma)$ over some finite time interval may be guaranteed by a large deviation principle \cite{FW98}. 
However, large deviation principles only apply in the small noise regime. 

When $\tau_\sigma$ is almost surely finite, but large, we might study the dynamics of the system in terms of a ``quasi-ergodic measure''. 
The textbook of Collet \emph{et al.}~\cite{CMSM12}~or the recent dissertation of Villemonais \cite{V19}~are good entry points to the general theory of quasi-ergodic measures. 
An extensive bibliography of works related to quasi-ergodic measures has been collected by Pollett \cite{P15}. \footnote{When introducing quasi-ergodic measures, it is typical -- and important -- to distinguish between quasi-ergodic measures and quasi-stationary measures. 
However, we choose to omit a discussion of these and other subtleties, opting to only state immediately necessary definitions and results. } 

\begin{Def}
Let $(E,\mathcal{E})$ be a metric space, $(Y_t)_{t\ge0}$ an $E$-valued Markov process, and $B\subset E$ a bounded open subset of $E$. 
Suppose $Y_0\in B$, and define the exit time of $(Y_t)_{t\ge0}$ as 
\[
\tau_B\,\coloneqq\,\inf\{t>0\,:\,Y_t\in\partial B\}. 
\]
A \emph{quasi-ergodic measure}~of $(Y_t)_{t\ge0}$ in $B$ is a measure $\mu$ such that  
\[
\lim_{t\rightarrow\infty}\mathbb{E}_{\zeta}\left[\frac{1}{t}\int_0^t\chi_A(Y_s)\,ds\,|\,\tau_B>t\right]\,=\,\mu(A)\quad\forall A\in\mathcal{B}(E), 
\]
for any initial distribution $\zeta$ supported in the open set $B$. 
That is, $\mu(A)$ is the expected  fraction of time $(Y_t)_{t\ge0}$ spends in $A$, given that $Y_t\in B$ for all $t\ge0$. 
Note that a quasi-ergodic measure is an ergodic measure in the usual sense if and only if $\tau_B=\infty$ almost surely. 
\end{Def}

To prove an analogue of Theorem \ref{thm:ErgodicFrequency}~in the quasi-ergodic setting, we need a version of the pointwise ergodic theorem for quasi-ergodic measures. 
The earliest ``quasi-ergodic theorem'', to our knowledge, was proven by Breyer \&~Roberts in 1999 \cite{BB99}. 
We make use of the quasi-ergodic theorem in \cite{V19}, modified from \cite{ZLS13}. 
Its statement is here modified to fit our notation. 


\begin{prop}[Example 5.1 and Corollary 6.5 of \cite{V19}]
\label{prop:quasiErgodic}
Let $E$ be a Banach space. 
Take a bounded subset $B$ of $E$, and let $(Y_t)_{t\ge0}$ be an $E$-valued It{\^o}~diffusion with locally H{\"o}lder continuous drift \&~diffusion coefficients, and a locally uniformly elliptic diffusion coefficient. 
Then, $(Y_t)_{t\ge0}$ possesses a unique quasi-ergodic distribution $\nu$. 
Moreover, for any $\nu$-measurable function $g:E\rightarrow\R$, any starting point $x\in B$, and any $\epsilon>0$, 
\begin{equation}
\label{eq:QuasiErgodic}
\pr_x\left(\left|\frac{1}{t}\int_0^tg(X_s)\,ds - \int_{B(\Gamma)}g(x)\,\nu(dx)\right|\ge\epsilon\,\,\big|\,\,t<\tau_\sigma\right)\,\xrightarrow[t\rightarrow\infty]{}\,0.
\end{equation}
\end{prop}

From Proposition \ref{prop:quasiErgodic}, we have the following.  

\begin{thm}
\label{thm:QuasiFrequency}
Let $(X_t^\sigma)_{t\ge0}$ be a stochastic process governed by \eqref{eq:SDE}. 
Let $\Gamma$ be a stable limit cycle of \eqref{eq:ODE}~with basin of attraction $B(\Gamma)$ and period $T>0$. 
Fix a phase map $\pi:B(\Gamma)\rightarrow[0,T)$. 
Suppose that 
\begin{enumerate}[(i)]
\item$(X_t^\sigma)_{t\ge0}$ satisfies the hypotheses of Proposition \ref{prop:quasiErgodic}, so that it has a quasi-ergodic measure $\mu_\sigma$ with bounded support in $B(\Gamma)$, 
\item If $\pi$ has a phase singularity $x_0\in B(\Gamma)$, then $\mu_\sigma$ is such that for some $\delta>0$ 
\[
\int_{B_\delta(x_0)}\norm*{\pi'(x)}\,\mu_\sigma(dx)\,<\,\infty,\qquad \int_{B_\delta(x_0)}\norm*{\pi''(x)}\,\mu_\sigma(dx)\,<\,\infty. 
\]
\end{enumerate}
Then, 
\begin{equation}
\label{eq:QuasiConvergence}
\pr_x\left(\left|\frac{1}{t}\pi_1(X_t^\sigma) - c_\sigma \right|\ge\epsilon\,\big|\,t<\tau_\sigma\right)\,\xrightarrow[t\rightarrow\infty]{}\,0, 
\end{equation} 
where 
\begin{equation}
\label{eq:QuasiFrequency} 
c_\sigma\,=\,\int_{B(\Gamma)}\pi_1'(x)V(x) + \frac{1}{2}\Tr\pi_1''(x)[B(x),B(x)]\,\mu_\sigma(dx). 
\end{equation}
\begin{proof}
For this proof, we fix $\sigma>0$ and write $X_t=X_t^\sigma$. 
Recall the It{\^o}~formula \eqref{eq:piIto}~for $\pi(X_t)$ and the local martingales $(I_t)_{t\ge0},\,(II_t)_{t\ge0}$, defined therein. 
The proof is similar to that of Theorem \ref{thm:ErgodicFrequency}, with some subtlety arising in the handling of $II_t$. 

When $(X_t)_{t\ge0}$ satisfies the hypotheses of Proposition \ref{prop:quasiErgodic}, we know that it has a well-defined $Q$-process in $B(\Gamma)$, which we denote $\tilde{X}_t$ (see Chapter 6 of \cite{V19}). 
We then define 
\[
\tilde{II}_t\,\coloneqq\,\int_0^t\pi'(\tilde{X}_s)B(\tilde{X}_s)\,dW_s. 
\]
As this is a local martingale with respect to $(\mathbb{Q}_x)_{x\in B(\Gamma)}$, we may apply It{\^o}'s isometry and the Burkholder-Davis-Gundy inequality to it. 
In particular, if the initial distribution of $\tilde{X}_t$ is $\mu_\sigma$, then 
\[
\begin{aligned}
c\,\mathbb{E}_{\mu_\sigma}\left[[\tilde{II}]_t\right]\,&\le\,\mathbb{E}_{\mu_\sigma}\left[\norm*{\tilde{II}_t}^2\right]\\
&\le\,\int_\Omega\int_0^t\norm*{\pi'(\tilde{X}_s(\omega))B(\tilde{X}_s(\omega))}^2\,ds\,\mathbb{Q}(d\omega)\\
&\le\, \int_0^t\,ds\,\int_{B(\Gamma)}\norm*{\pi'(x)B(x)}\,\mu_\sigma(dx) \,\le\, C_0\,t 
\end{aligned}
\]
for some $c,\,C_0>0$. 
A similar estimate can be obtained for an arbitrary initial distribution $\nu$ with support in $B(\Gamma)$, since the distribution of $\tilde{X}_t$ becomes arbitrarily close to $\mu_\sigma$ for large $t>0$. 
Applying Theorem 4.1 of van Zanten \cite{vZ00}, we conclude that, 
for any $\epsilon>0$ and initial distribution $\nu$ supported in $B(\Gamma)$, 
\begin{equation}
\label{eq:QII}
\mathbb{Q}_\nu\left(\frac{1}{t}\tilde{II}_t>\epsilon\right)\,\xrightarrow[t\rightarrow\infty]{}\,0. 
\end{equation}
By the definition of $\mathbb{Q}_\nu$ (see Chapter 6 of \cite{V19}), for any $\epsilon_0>0$ and $A\in\mathcal{F}$ there exists $t_0>0$ such that if $t>t_0$, then 
\[
\norm*{\mathbb{Q}_\nu(A) - \pr_\nu(A\,|\,t<\tau_\sigma)}\,<\,\epsilon_0. 
\]
Hence, \eqref{eq:QII}~implies 
\[
\pr_\nu\left(\frac{1}{t}II_t>\epsilon\,|\,t<\tau_\sigma\right)\,=\,\pr_\nu\left(\frac{1}{t}\tilde{II}_t>\epsilon\,|\,t<\tau_\sigma\right)\,\xrightarrow[t\rightarrow\infty]{}\,0. 
\]
Applying Theorem \ref{prop:quasiErgodic}~to $(I_t)_{t\ge0}$ with respect to the quasi-ergodic measure $\mu_\sigma$ completes the proof. 
\end{proof}
\end{thm}

\begin{Def}
If \eqref{eq:SDE}~satisfies the hypotheses of Theorem \ref{thm:QuasiFrequency}, we refer to $c_\sigma$ defined in \eqref{eq:QuasiFrequency}~as the \emph{quasi-asymptotic frequency}~of \eqref{eq:SDE}. 
\end{Def}

Though its proof is simple, Theorem \ref{thm:QuasiFrequency}~is stronger than the results of \emph{e.g.}~Giacomin \emph{et al.}~\cite{GPS18}, in that we are not restricted to the small $\sigma$ regime. 
Even in the small $\sigma$ regime, we will see in Section \ref{sec:Comparison}~that our result may differ from that of \cite{GPS18}. 
To bridge this (potential) gap between our prediction and \cite{GPS18}, we remark that when the quasi-asymptotic frequency does exist, it may only be observable on a particular time scale. 
We therefore expect that, when our prediction for $c_\sigma$ differs from the prediction of Giacomin \emph{et al.}, the two predictions are observable on different time scales. 

It should also be noted that, even when the quasi-ergodic frequency in \eqref{eq:QuasiFrequency}~exists, it is not necessarily observable on any time scale with high probability. 
More precisely, the rate of convergence of $t^{-1}\pi_1(X_t^\sigma)$ to $c_\sigma$ may be slower than the rate of escape of $X_t^\sigma$ from $B(\Gamma)$. 
In this case $c_\sigma$ has a low probability of being observed. 

Unfortunately, the theory of quasi-ergodic measures is still insufficiently developed for a general comparison of the rate of convergence of $t^{-1}\pi_1(X_t^\sigma)$ with the rate of escape of $X_t^\sigma$ from $B(\Gamma)$. 
We can nevertheless find several examples where Theorem \ref{thm:QuasiFrequency}~applies, and where the rate of convergence to the quasi-asymptotic frequency appears to be much faster than the rate of escape from $B(\Gamma)$. 

\begin{figure}[h!]
\centering
\includegraphics[width=0.5\textwidth]{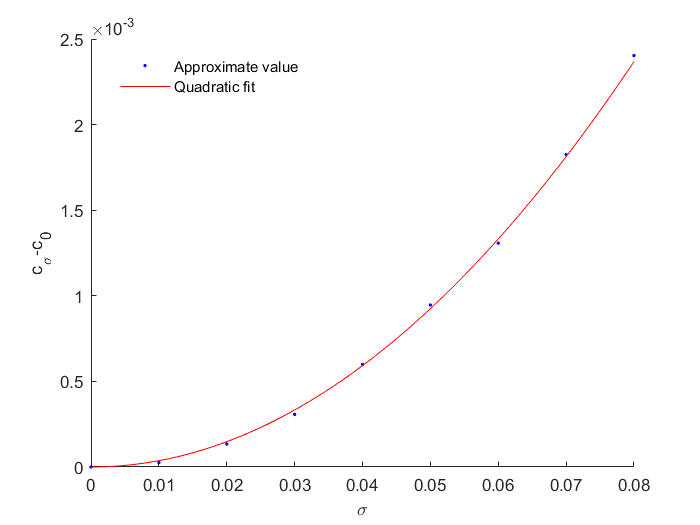}
\caption{ 
(a) Approximate values of the difference between the quasi-asymptotic frequency $c_\sigma$ and the deterministic frequency $c_0$ in \eqref{eq:3Cycles}. 
Numerical parameters were $a=1,\,b=2,\,c=3$. 
The system was simulated up to time $t_{end}=400000$ using an Euler-Maruyama scheme with time step $dt=0.0004$. 
Simulating on this time scale, $t^{-1}\pi_1(X_t)$ appears to converge to the quasi-asymptotic frequency, and the approximation $c_\sigma-c_0\simeq m\sigma^2$ appears to be good for $m=0.37$. 
} 
\label{fig:3Cycles}
\end{figure}

\begin{EX}
\label{EX:3Cycles}
Consider the following SDE in polar coordinates, with deterministic parameters $0<a<b<c$, 
\begin{equation}
\label{eq:3Cycles} 
\begin{aligned}
dr\,&=\,r(r-a)(r-b)(r-c)\,dt + \sigma r\,dW, \\
d\theta\,&=\,r^{-2}\,dt. 
\end{aligned}
\end{equation}
The system with $\sigma=0$ has a stable limit cycle at $r=b$, the basin of attraction of which is bounded by unstable limit cycles at $r=a$ and $r=c$. 
Example 5.1 of \cite{V19}~guarantees the existence of a unique quasi-ergodic measure supported in the annulus $A_{a,c}\coloneqq\{(r,\theta)\,:\,r\in(a,c)\}$. 
There is no phase singularity on the boundary of this domain. 
Hence, \eqref{eq:3Cycles}~possesses a quasi-asymptotic frequency in $A_{a,c}$ for $\sigma>0$. 

We have numerically approximated the quasi-ergodic frequency of \eqref{eq:3Cycles}, and show the results in Figure \ref{fig:3Cycles}. 
For $\sigma\in[0,0.08]$, \eqref{eq:3Cycles}~was simulated up to time $t_{end}=400000$ using an Euler-Maruyama scheme (see Chapter 8 of Lord \emph{et al.}~\cite{LPS14}) with step size $dt=0.0004$. 
Note that $c_\sigma$ appears to depend quadratically on $\sigma$ well into the moderate-noise regime (see Figure \ref{fig:3Cycles}). 
Convergence of the quasi-asymptotic frequency appears to occur slowly, but still much faster than escape from $A_{a,c}$. 
We were not able to simulate long enough to observe a sufficient number of escapes to estimate the escape rate. 
\end{EX}

\begin{figure}[h!]
\centering
\includegraphics[width=0.5\textwidth]{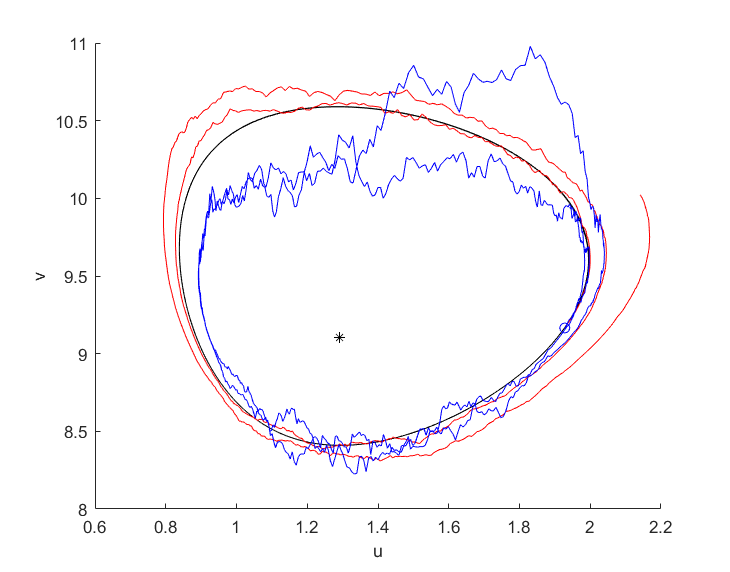}
\caption{ 
Trajectories of the predator-prey system \eqref{eq:PP}~with parameters as in \eqref{eq:PPparam}. 
In black: the stable limit cycle of the system when $\sigma=0$ surrounding the ustable equilibrium $(u_*,v_*)$, the latter shown as a star. 
In red: a trajectory of the system started on the limit cycle with $\sigma=0.1$ and $B=B_0$. 
In blue: a trajectory of the system started on the limit cycle with $\sigma=0.5$ and $B=B_0$. 
} 
\label{fig:PPphase}
\end{figure}

\begin{figure}[h!]
\centering
\begin{subfigure}[b]{0.45\textwidth}
\centering
\includegraphics[width=\textwidth]{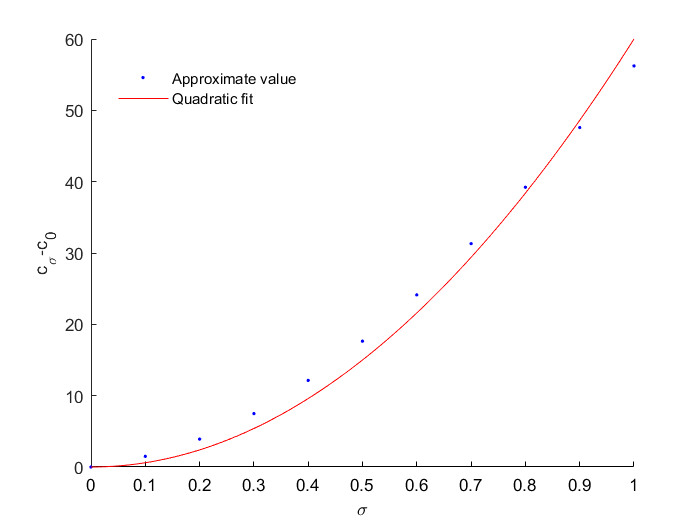}
\caption{}
\label{fig:PP_B0_frequencies}
\end{subfigure}
\begin{subfigure}[b]{0.45\textwidth}
\centering
\includegraphics[width=\textwidth]{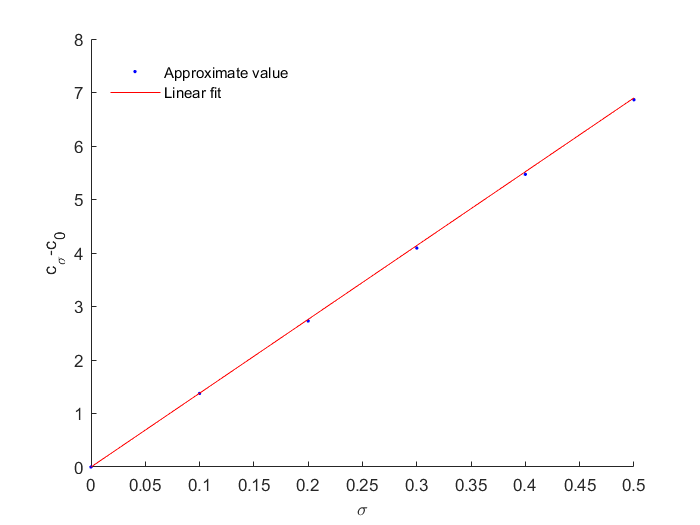}
\caption{}
\label{fig:PP_B1_frequencies}
\end{subfigure}
\caption{ 
Approximate values of $c_\sigma$ for $\sigma\in\{0,0.1,\ldots,1\}$ in \eqref{eq:PP}, with parameters as in \eqref{eq:PPparam}. 
In (a), $B=B_0$. 
In (b), $B=B_1$. 
All simulations were performed with a time step size of $dt=0.01$, up to time $t_{end}=1000000$. 
Note that for $B=B_0$, $c_\sigma$ does not exhibit a quadratic dependence on $\sigma>0$ over this range of values. 
Meanwhile, for $B=B_1$, $c_\sigma$ appears to depend linearly on $\sigma>0$ over the range of values shown. 
However, see Remark \ref{rmk:Numerics}. 
} 
\label{fig:PPcsigma}
\end{figure}

\begin{figure}[h!]
\centering
\begin{subfigure}[b]{0.45\textwidth}
\centering
\includegraphics[width=\textwidth]{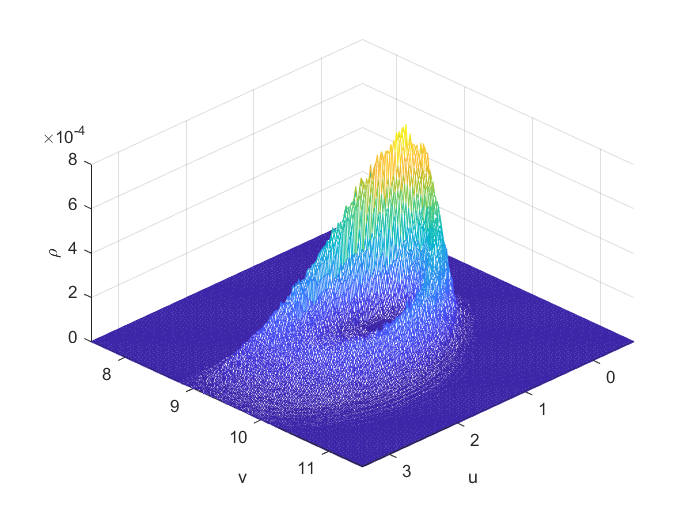}
\caption{}
\label{fig:PP_B0_SmallNoise_measure}
\end{subfigure}
\begin{subfigure}[b]{0.45\textwidth}
\centering
\includegraphics[width=\textwidth]{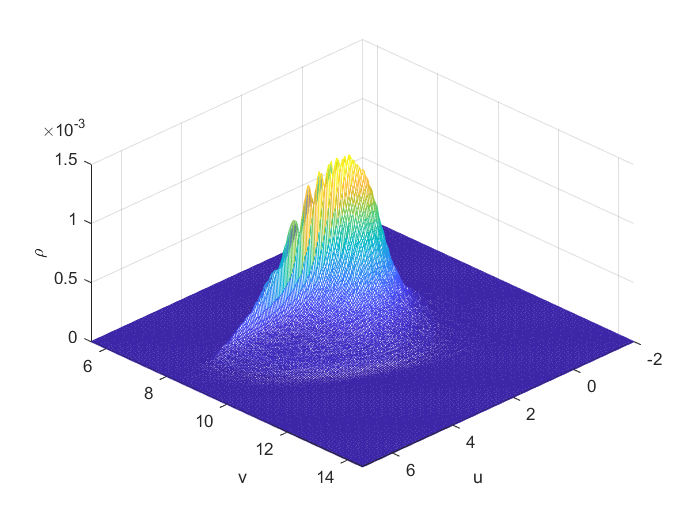}
\caption{}
\label{fig:PP_B0_LargeNoise_measure}
\end{subfigure}
\begin{subfigure}[b]{0.45\textwidth}
\centering
\includegraphics[width=\textwidth]{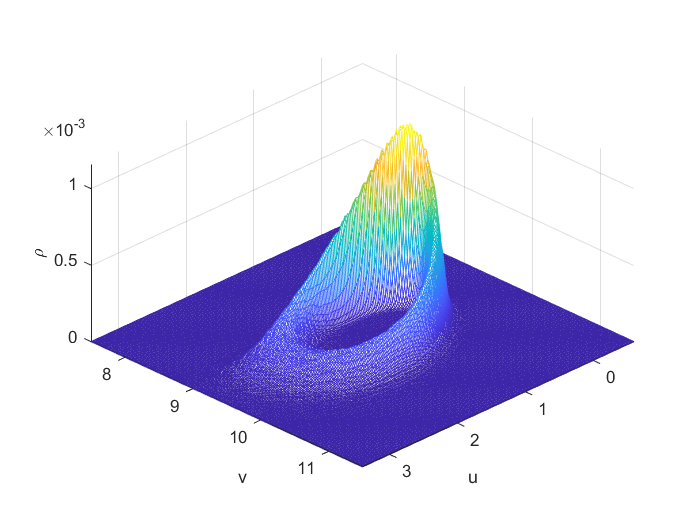}
\caption{}
\label{fig:PP_B1_SmallNoise_measure}
\end{subfigure}
\begin{subfigure}[b]{0.45\textwidth}
\centering
\includegraphics[width=\textwidth]{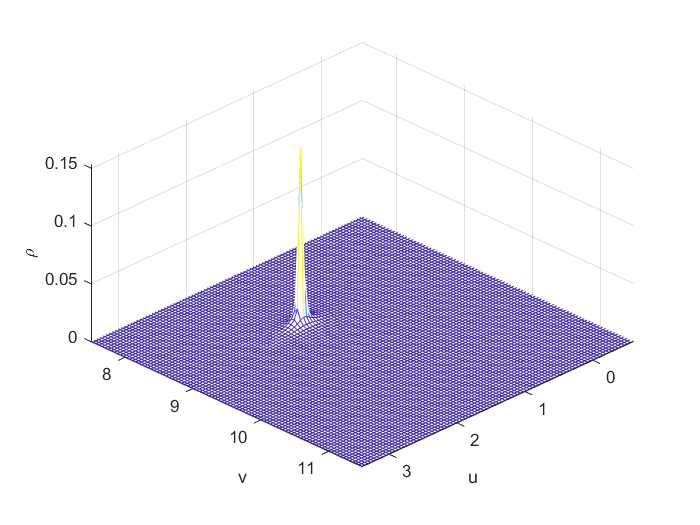}
\caption{}
\label{fig:PP_B1_LargeNoise_measure}
\end{subfigure}
\caption{ 
Monte-Carlo approximations of the quasi-stationary measures of \eqref{eq:PP}, with parameters as in \eqref{eq:PPparam}. 
The diffusion coefficient in each subfigure is as follows: 
(a) $\sigma=0.1$, $B=B_0$; 
(b) $\sigma=0.5$, $B=B_0$; 
(c) $\sigma=0.1$, $B=B_1$; 
(d) $\sigma=0.5$, $B=B_1$. 
The introduction of noise on this scale has a dramatic effect on the quasi-stationary measure of the system. 
In \ref{fig:PP_B0_LargeNoise_measure}, the apparent multiple peaks are a numerical artefact. 
In \ref{fig:PP_B1_SmallNoise_measure}~\&~\ref{fig:PP_B1_LargeNoise_measure}, we see that the choice of diffusion coefficient keeps the system, on average, further away from the phase singularity at $(u_*,v_*)$. 
In \ref{fig:PP_B1_LargeNoise_measure}, the system spends most of its time near the peak shown in the figure. 
However, by simulating individual paths of the system, we see that at random times \eqref{eq:PP}~still makes rapid circuits around $(u_*,v_*)$. 
Thus the system is still oscillatory, though this is not apparent from the form of its invariant measure when approximated at the coarse scale used here. 
} 
\label{fig:PPmeasures}
\end{figure}

\begin{EX}
\label{EX:PP} 
We here consider the example of a stochastic predator-prey system with a Holling type III functional response, 
\begin{equation}
\label{eq:PP}
\begin{aligned}
du\,&=\, \left(u(a-u) - b\frac{u^2v}{1+u^2}\right)\,dt \\
dv\,&=\, \left(c\frac{u^2v}{1+u^2} - dv\right)\,dt + \sigma B(v)\,dW. 
\end{aligned}
\end{equation}
In applications, \eqref{eq:PP}~has been used to model the interaction of a population of phytoplankton, $u$, with a population of zooplankton, $v$. 
The parameters of \eqref{eq:PP}~are non-dimensionalized, but derive from the birth, death, and feeding rates of both populations, along with the carrying capacity of their environment. 
We refer to Freedman \cite{F80}~for a further discussion of the deterministic version of this model.  
For the stochatic model see, for instance, Reichenbach \emph{et al.}~\cite{RMF08}~or Sun \emph{et al.}~\cite{S10}. 
We remark that systems such as \eqref{eq:PP}~are often studied with spatial diffusion and a time-periodic forcing term. 
For simplicity, we omit these features. 

We simulate \eqref{eq:PP}~with two different choices of diffusion coefficient, taking $B$ to be either 
\[
B_0(v)\,\coloneqq\,1\qquad\text{ or }\qquad B_1(v)\,\coloneqq\,v-v_*. 
\]
The noise $B=B_0$ is motivated by the analysis of Reichenbach \emph{et al.}~\cite{RMF08}, while we also study the noise $B=B_1$ out of curiousity. 
Either choice of noise may cause the system to exit the positive quadrant of $\R^2$ in finite time, beyond which point the model ceases to be meaningful. 
The model \eqref{eq:PP}~is usually analyzed outside of the small-noise regime, so that the results of those cited in Section \ref{subsec:Background}~cannot apply. 
Hence, we use the theory of quasi-ergodic mesaures, conditioning on the system remaining in the positive quadrant of $\R^2$. 

The deterministic parameter values we use are 
\begin{equation}
\label{eq:PPparam} 
\begin{aligned}
a\,=\, 6.8,\quad b\,=\,1.25,\quad c\,=\,0.8,\quad d\,=\,0.5. 
\end{aligned}
\end{equation}
At these parameter values, the system with $\sigma=0$ possesses an unstable equilibrium at 
\[
u_*\,\coloneqq\,\frac{\sqrt{d(c-d)}}{c-d}\,\simeq\,1.29,\qquad v_*\,\coloneqq\,\frac{c(a-u_*)}{bu_*(c-d)}\,\simeq\,9.10,
\]
surrounded by a stable limit cycle $\Gamma$. 
See Figure \ref{fig:PPphase}. 
Since the isochronal phase of the stable limit cycle in \eqref{eq:PP}~is not easy to compute numerically, we use the angular phase centered at $(u_*,v_*)$; that is, we take 
\[
\pi(u,v)\,\coloneqq\,\tan\left(\frac{u-u_*}{v-v_*}\right). 
\]
In the case of \eqref{eq:PP}, the angular phase centered at $(u_*,v_*)$ satisfies our definition of a phase map. 
Note that $(u_*,v_*)$ is a phase singularity of $\pi$. 

In this setup, all hypotheses of Theorem \ref{thm:QuasiFrequency}~are automatically satisfied, save for (ii). 
We have not been able to work out sufficient conditions for (ii) to be satisfied, such as a quasi-ergodic analogue of Theorem \ref{cor:}~(but see the discussion at the end of Appendix \ref{app:Sufficient}). 
Nevertheless, we may numerically observe an apparent convergence of $t^{-1}\pi_1(X_t^\sigma)$ to a deterministic constant for a range of values of $\sigma>0$ and either choice of $B$, as we discuss in the following. 

When $B=B_0$, we see that the response of $c_\sigma$ to increasing $\sigma>0$ is large, and apparently non-quadratic when observed on this scale. 
See Figure \ref{fig:PP_B0_frequencies}. 
The noise has a correspondingly large qualitative effect on the quasi-stationary measure of the system (and thus the closely related quasi-ergodic measure) -- see Figures \ref{fig:PP_B0_SmallNoise_measure}~\&~\ref{fig:PP_B0_LargeNoise_measure}. 
Escapes from the positive quadrant of $\R^2$ were not observed on the time-scales which we simulated over. 

When $B=B_1$, increasing $\sigma$ causes a large increase in $c_\sigma$ until about $\sigma=0.5$. 
The dependence of $c_\sigma$ on $\sigma$, viewed at this scale, is apparently linear, with a slope of $m\simeq13.8$. 
For values of $\sigma$ greater than $0.5$, the choice of diffusion coefficient caused the system to escape from the positive quadrant of $\R^2$ before suitable convergence was observed.  
In this situation, available computer time would not allow us to approximate $c_\sigma$ using the methods of this paper. 
Simultaneously, the choice $B=B_1$ kept the system, on average, further away from  the phase singularity at $(u_*,v_*)$ than with $B=B_0$. 
See Figure \ref{fig:PPmeasures}. 
This may account for the smaller effect which increasing $\sigma$ had on the quasi-asymptotic frequency relative to when $B=B_0$, as we see in Fgure \ref{fig:PPcsigma}. 

For both $B=B_0$ and $B=B_1$, the rate of convergence to $c_\sigma$ was much slower in this example than any of the other examples we have discussed. 
Nevertheless, with $B=B_0$ and $\sigma<1$, escapes from the positive quadrant of $\R^2$ were not observed, so that $c_\sigma$ can be approximated by $\frac{1}{t}\pi_1(X_{t})$ for large values of $t$. 
The same holds for $B=B_1$ when $\sigma<0.5$. 
For $B=B_1$ and $\sigma>0.5$, the rate of escape from the positive quadrant of $\R^2$ was sufficiently fast to prevent us from observing $\frac{1}{t}\pi_1(X_t)$ for sufficiently large $t$, due to limited computer time. 
\end{EX}

\begin{rmk}
\label{rmk:Numerics}
Our approach to computing $c_\sigma$ is not suitable for studying the small noise regime. 
This is due to the fact that, if we want an accurate simulation of an SDE, we need to take the time step $dt$ at least on the order of $\sigma^2$. 
Meanwhile, if we want an accurate approximation of $c_\sigma$, we need to simulate for a very long time. 
Hence, our approach to approximating $c_\sigma$ for small $\sigma>0$ would require more computer time than we have at our disposal. 
Future work may attempt to compute $c_\sigma$ directly from \eqref{eq:ErgodicFrequency}~by numerically approximating the quasi-ergodic measure $\mu_\sigma$, and then numerically integrating. 

Unfortunately, numerical techniques for approximating quasi-ergodic measures are still in their infancy, but see Dobson \emph{et al.}~\cite{DLZ19}, Li~\cite{L18}, or Li \&~Yuan \cite{LY21}. 
The Monte-Carlo method, which we have used in our approximations of stationary and quasi-stationary measures, does not provide sufficient accuracy for an accurate computation of $c_\sigma$. 
Our approximations are only useful insofar as they provide a qualitative understanding of the co-responsive behaviour of $c_\sigma$ and $\mu_\sigma$ to changes in $\sigma$. 
\end{rmk}

\begin{rmk}
\label{EX:FHNa} 
The theory of quasi-ergodic measures has been developed in the setting of Markov processes taking values in an arbitrary measurable space \cite{V19}.
Hence, just as in Example \ref{EX:FHN}, the theory of this section extends to the infinite dimensional setting, with the added difficulties associated with invariant measure theory in spaces which are not locally compact. 
We could for instance consider the FitzHugh-Nagumo system perturbed by additive noise, 
\begin{equation}
\label{eq:FHNa}
\begin{aligned}
du\,&=\, \big(D\Delta u + u(1-u)(u-a)\big)\,dt + \sigma\,dW, \\
dv\,&=\, \big(\delta D\Delta v + \gamma u - v\big)\,dt, 
\end{aligned}
\end{equation}
again considered as an evolution equation on $L^2(S^1;\R^2)$. 

Just as for ergodic measures, proving the existence of a quasi-ergodic measure in the infinite dimensional setting is somewhat more subtle than in the finite dimensional case. 
At this point, there are no easily applied theorems in the quasi-ergodic setting which would be analogous to Theorem 11.38 of Da Prato \&~Zabczyk \cite{DPZ14}. 
One could attempt to restrict dynamics to a \emph{compact}~subdomain of phase space, but the escape rate of a truly infinite dimensional process from such a domain is likely to be high. 
Further study is needed. 
\end{rmk}

\section{Decomposing the asymptotic frequency \&~comparison with past results}
\label{sec:Comparison} 
We have seen that if \eqref{eq:SDE}~has a unique (quasi) ergodic measure $\mu_\sigma$ in $B(\Gamma)$, then it has a (quasi) asymptotic frequency in $B(\Gamma)$. 
For a fixed phase map $\pi$, this frequency is given by 
\begin{equation}
\label{eq:ErgodicFrequency2}
c_\sigma\,=\,\int_{B(\Gamma)}\pi_1'(x)V(x) + \frac{\sigma^2}{2}\Tr\pi_1''(x)[B(x),B(x)]\,\mu_\sigma(dx), 
\end{equation} 
so long as the integral converges. 
In this section, we discuss the possibility of a non-quadratic response of $c_\sigma$ to $\sigma$. 
Additionally, we compare \eqref{eq:ErgodicFrequency2}~with a prediction of Giacomin \emph{et al.}~\cite{GPS18}. 

Remark that for any continuous function $f$ which is uniformly bounded on a neighbourhood of $\Gamma$, we have 
\[
\int_{B(\Gamma)} f(x)\,\mu_0(dx)\,=\,\frac{1}{T}\int_0^T f(\gamma_t)\,dt, 
\]
where we recall that $\{\gamma_t\}_{t\in\R}$ is the parameterization of $\Gamma$ of period $T>0$. 
Decomposing the (quasi) ergodic measure of \eqref{eq:SDE}~in $B(\Gamma)$ as $\mu_\sigma = \mu_0+\nu_\sigma$ for some signed measure $\nu_\sigma$, we write the (quasi) asymptotic frequency of \eqref{eq:SDE}~-- assuming it exists -- as 
\begin{equation}
\label{eq:AsymptoticFrequencyDecomposition} 
\begin{aligned}
c_\sigma\,&=\,\frac{1}{T}\int_0^T \pi_1'(\gamma_t)V(\gamma_t)\,dt + \frac{\sigma^2}{2T}\int_0^T\Tr\pi_1''(\gamma_t)[B(\gamma_t),B(\gamma_t)]\,dt \\
&\qquad\, +\int_{B(\Gamma)} \pi'(x)V(x)\,\nu_\sigma(dx) + \frac{\sigma^2}{2}\int_{B(\Gamma)}\Tr\pi_1''(x)[B(x),B(x)]\,\nu_\sigma(dx) \\
&\eqqcolon\, c_0 + a_\sigma + \sigma^2(b_0  + b_\sigma).  
\end{aligned}
\end{equation}

When $\pi$ is the isochron map, it holds that $\pi'(x)V(x)=1$ for $x\in B(\Gamma)$, so 
\[
\int_0^T\pi_1'(\gamma_t)V(\gamma_t)\,dt \,=\,1\quad \text{ and }\quad
\int_{B(\Gamma)}\pi_1'(x)V(x)\,\nu_\sigma(dx)\,=\,0. 
\]
Hence, when $\pi$ is the isochron map \eqref{eq:AsymptoticFrequencyDecomposition}~becomes  
\begin{equation}
\label{eq:IsochronalFrequencyDecomposition} 
\begin{aligned}
c_\sigma\,&=\, 
c_0 + \sigma^2(b_0  + b_\sigma). 
\end{aligned}
\end{equation}
A quadratic approximation of $(c_\sigma)_{\sigma\ge0}$ is therefore good if $b_\sigma$ remains relatively constant. 
This would follow, for instance, from a bound on the total variation norm of $\nu_\sigma$. 

We can compare \eqref{eq:IsochronalFrequencyDecomposition}~with a prediction of \cite{GPS18}, who study the long-time behaviour of the isochronal phase of SDE of the form \eqref{eq:SDE}. 
Notably, \cite{GPS18}~do not require $\tau_\sigma=\infty$, and handle the fact that the isochronal phase may not be defined for all time by taking $\sigma\rightarrow0$ simultaneously with $t\rightarrow\infty$. 
That is, the authors take a sequence of times $(t(\sigma))_{\sigma>0}$ such that 
\[
t(\sigma)\xrightarrow[\sigma\rightarrow0]{}\infty  
\]
in a controlled fashion, and study the behaviour of $X_{t(\sigma)}^\sigma$ as $\sigma\rightarrow0$. 
The family $(t(\sigma))_{\sigma>0}$ is chosen such that $X_{t(\sigma)}^\sigma$ remains in $B(\Gamma)$ with high probability for each $\sigma>0$, thanks to large deviation estimates. 

To state a result of \cite{GPS18}, we recall the \emph{winding number}~of $X_t^\sigma$, which is the number of full clockwise rotations minus the number of full counter-clockwise rotations made by $\pi(X_t^\sigma)$. 
The winding number at time $t\ge0$ is denoted $n_t^\sigma$. 

\begin{thm}[Theorem 2.6 of \cite{GPS18}]
\label{thm:GPSlarge}
There exists $c>0$ such that for any $(t(\sigma))_{\sigma\ge0}$ satsifying 
\[
\lim_{\sigma\rightarrow0+}\sigma^2t(\sigma)\,=\,\infty,\qquad \lim_{\sigma\rightarrow0+}e^{-c\sigma^{-2}}t(\sigma)\,=\,0, 
\]
then 
\begin{equation}
\label{eq:ResultGPS}
\lim_{\sigma\rightarrow0+}\frac{n^\sigma_{t(\sigma)}-t(\sigma)/T}{\sigma^2t(\sigma)/T}\,=\,b_0,
\end{equation}
with $b_0$ as in \eqref{eq:AsymptoticFrequencyDecomposition}.  
\end{thm}

Colloquially, Theorem \ref{thm:GPSlarge}~implies that for ``small'' $\sigma>0$ and ``appropriately sized'' $t>0$ (specifically, such that $\sigma^{-2}\lesssim t\lesssim e^{-c\sigma^{-2}}$), the time average frequency of the random oscillator at time $t$, denoted $\tilde{c}_\sigma(t)$, is 
\[
\tilde{c}_\sigma(t)\,\coloneqq\,\frac{n_t^\sigma}{t}\,\simeq\,c_0 + \sigma^2b_0. 
\]
That is, the stochastic frequency is, approximately in the long term, equal to the deterministic frequency plus a correction of order $\sigma^2$. 
Using a zero-one law, it can be shown that $n_t^\sigma/t\simeq \pi_1(X_t^\sigma)/t$ for large $t\ge0$ (assuming $t\le\tau_\sigma$). 
Hence, we find that our result agrees with the result of \cite{GPS18}~in the small $\sigma>0$ regime, unless $b_\sigma$ in \eqref{eq:AsymptoticFrequencyDecomposition}~varies significantly with increasing $\sigma$. 

It would be interesting to find an example of \eqref{eq:SDE}~where $b_\sigma$ changes rapidly for small $\sigma>0$, 
\[
\frac{\partial}{\partial\sigma}\big\rvert_{\sigma=0}b_\sigma\,\neq\,0. 
\]
In this case, our definition of the asymptotic frequency of a stochastic oscillator may differ from that of \cite{GPS18}~in the small noise regime. 
In such a scenario, it is likely that our $c_\sigma$ is observable on a different time scale than the $\tilde{c}_\sigma(t)$ predicted by \cite{GPS18}. 
Unfortunately, as we have already noted, the theory of quasi-ergodic measures remains insufficiently devloped for us to estimate the time-scale on which $c_\sigma$ is observed with high probability. 

Further understanding of the relation between $c_\sigma$ and $\mu_\sigma$ can only be achieved on a case-by-case basis, and requires more information about $\mu_\sigma$. 
All simulations used in this paper are the result of Monte Carlo schemes.  
In the future, improved methods of approximating the invariant measures of stochastic oscillators could be used to obtain more accurate predictions of $c_\sigma$ using the formula \eqref{eq:ErgodicFrequency}, such as Dobson \emph{et al.}~\cite{DLZ19}~or Li~\cite{L18}. 
A similar technique has been developed for approximating quasi-ergodic measures in Li \&~Yuan \cite{LY21}.

\section{Conclusions \&~outlook}
\label{sec:Conclusions} 
The response of the asymptotic frequency of a stochastic oscillator to varying noise amplitude has been studied extensively over the past two decades. 
Before defining the asymptotic frequency of a stochastic oscillator, one needs a notion of its phase. 
This is usually achieved via a so-called phase map, as described in Section \ref{sec:Setup}. 
However, as the domain of definition of a phase map is often bounded, and stochastic oscillators driven by many realistic choices of noise will almost surely exit any bounded domain in finite time, the asymptotic frequency is not necessarily well-defined. 

When a stochastic oscillator remains in the domain of a phase map's definition for all time, one can obtain a formula for its asymptotic frequency via the pointwise ergodic theorem. 
This requires the phase map and ergodic measure of the stochastic oscillator to satisfy the integrability condition of Theorem \ref{thm:ErgodicFrequency}~(or Remark \ref{rmk:Unbounded}). 
While these conditions are satisfied in most applications, the possibility of pathological examples where they are not have not been ruled out. 

The theory of quasi-ergodic measures can be used to make sense of the asymptotic frequency of any stochastic oscillator. 
To avoid issues caused by the unboundedness of some stochastic processes, we may restrict our attention to the event of a stochastic oscillator remaining in a bounded subdomain of its phase space for all time. 
Though this event may be of probability zero, it can still yield physically relevant information. 
Whether or not it does depends on the comparative values of two quantities: 
the rate of convergence of the frequency of the (conditioned) stochastic oscillator to its asymptotic average, and the rate of escape of the stochastic oscillator from the bounded subdomain. 
If the latter is less than the former, then the quasi-asymptotic frequency is observable. 
Identifying conditions which guarantee the observability of quasi-asymptotic frequencies will be the subject of future work. 

Most previous studies of the response of the asymptotic frequency of stochastic oscillators to varying noise amplitude imply that the frequency depends quadratically on the noise amplitude, at least in the small noise regime. 
The results of Section \ref{sec:Comparison}~suggest that one could expect a quadratic response of the (quasi) asymptotic frequency both within and outside of the small noise regime, so long as the term $b_\sigma$ in \eqref{eq:AsymptoticFrequencyDecomposition}~remains relatively constant. 
On the other hand, it is possible that $(c_\sigma)_{\sigma>0}$ is not quadratic in the small noise regime if $b_\sigma$ varies significantly for small $\sigma$. 
This would require the difference between the (quasi) ergodic measure of the perturbed system to significantly differ from the ergodic measure of the unperturbed system. 
Determining when this is the case requires a better understanding of the dependence of the (quasi) ergodic ditribution of stochastic oscillators on noise amplitude, and this can only be achieved on a case-by-case basis.

\subsection*{Acknowledgements} 
Thanks to Maximilian Engel, for several helpful conversations during the preparation of this paper; 
to Rishabh Gvalani, for assistance with the proof of Proposition A.1;
to two anonymous referees, for their many helpful comments;  
and to J{\"u}rgen Jost, for his continuing patience and support. 
This work was funded by the International Max Planck Research School for Mathematics in the Sciences.

\appendix 
\section{Sufficient Conditions for Theorem \ref{thm:ErgodicFrequency}}
\label{app:Sufficient} 
The following proposition provides sufficient conditions for the hypotheses of Theorem \ref{thm:ErgodicFrequency}~to be satisfied. 
Due to the cursory nature of this paper, we have chosen to include only a sketch of its proof. 

\begin{prop}
\label{cor:}
Let $\pi$ be the isochron map of \eqref{eq:ODE}. 
Let $(X_t^\sigma)_{t\ge0}$ be the solution of \eqref{eq:SDE}. 
If $x_0\in B(\Gamma)$ is the only phase singularity of $\pi$ and  
\begin{enumerate}[(i)]
\item $B(x)$ and $V(x)$ have polynomial entries, and 
$\norm*{B(x)}$ and $\norm*{V(x)}$ are $O(\norm*{x-x_0})$ as $x$ approaches $x_0$, 
\item For all $t\ge0$, $X_t^\sigma$ has a density with respect to Lebesgue measure, 
\item $B(x)>0$ for $x\in B(\Gamma)$ and $B(x)\rightarrow0$ as $x$ approaches $\partial B(\Gamma)$, 
\item $\Tr V'(x_0)>2d\sigma^2$ and $V'(0)$ is positive definite, 
\end{enumerate}
then the assumptions of Theorem \ref{thm:ErgodicFrequency}~are satisfied. 
\begin{proof}
By assumption (iii) the process is trapped in $B(\Gamma)$ for all time, and we may therefore restrict our attention to the pre-compact Polish space $B(\Gamma)$.  
The existence of an ergodic measure of \eqref{eq:SDE}~in $B(\Gamma)$ then follows from a straightforward application of the Krylov-Bogoliubov Theorem and irreducibility. 
Its uniqueness is guaranteed by the irreducibility of the solution to \eqref{eq:SDE}~in $B(\Gamma)$, which also follows from assumption (iii). 
The ergodic measure has a density with respect to Lebesgue measure, denoted $\rho_\sigma$, by assumption (ii). 

Without loss of generality, we assume $x_0=0$. 
Since 
\begin{equation}
\label{eq:piV1}
\pi'(x)V(x)\,\equiv\,1\quad\text{ and }\quad \pi''(x)V(x) + \pi'(x)V'(x)\,\equiv\,0, 
\end{equation}
we have 
\[
\norm*{\pi'(x)}\,\sim\,\norm*{V(x)}^{-1}\quad\text{ and }\quad \norm*{\pi''(x)}\,\sim\,\norm*{\pi'(x)}\norm*{V'(x)}\norm*{V(x)}^{-1}. 
\]
By condition (i) of this corollary, this implies that, for $x$ near zero, 
\[
\norm*{\pi'(x)}\,\sim\norm*{x}^{-1}\quad\text{ and }\quad \norm*{\pi''(x)}\,\sim\,\norm*{x}^{-2}. 
\]
We will therefore aim to show that $\rho_\sigma(x)\le k\norm*{x}^2$ for some $k>0$. 

To this end, remark that the distribution $\rho_\sigma$ of $\mu_\sigma$ satisfies the stationary Fokker-Planck equation 
\begin{equation}
\label{eq:FP}
\nabla\cdot(V(x)\rho(x)) - \frac{\sigma^2}{2}\Delta(B(x)B(x)^*\rho(x))\,=\,0,\qquad\int_{B(\Gamma)}\rho_\sigma(x)\,dx\,=\,1. 
\end{equation}
Define $A\,\coloneqq\,2V'(x_0)/\sigma^2$. 
Fixing $\delta>0$ and taking arbitrarily small $\epsilon>0$, we use assumption (i) to approximate the solution $\rho_\sigma$ of \eqref{eq:FP}~by the solution $\tilde{\rho}_{\sigma,\epsilon}$ of the uniformly elliptic equation 
\begin{equation}
\label{eq:FPA}
\begin{aligned}
&\nabla\cdot(Ax\rho(x)) - \Delta((\norm*{x}^2+\epsilon)\rho(x))\,=\,0\quad\text{ for }\quad x\in B_\delta(0), \\
&\rho(x)\,=\,\rho_\sigma(x)\quad\text{ for }\quad x\in\partial B_\delta(0). 
\end{aligned}
\end{equation}
Since the solution $\rho_\sigma$ to \eqref{eq:FP}~is uniformly bounded on $\partial B_\delta(0)$, we can find $k>0$ such that 
\[
\rho_\sigma(x)\,\le\,k\norm*{x}^2\quad \text{ for }\quad x\in\partial B_\delta(0). 
\]
Writing \eqref{eq:FPA}~in divergence free form, a comparison principle (Theorem 3.3 of Gilbarg \&~Trudinger \cite{GT98}) applies if $\Tr A\ge 4d$ and $A$ is positive definite. 
This allows us to conclude that 
\[
\tilde{\rho}_{\sigma,\epsilon}(x)\,\le\,k\norm*{x}^2\quad\text{ for all }\quad x\in B_\delta(0), 
\]
where $k$ is independent of $\epsilon>0$. 

It can be shown that, for small enough $\epsilon>0$, $\tilde{\rho}_{\sigma,\epsilon}$ is indeed a good approximation of $\rho_\sigma$ on $B_\delta(x)$ when $\delta$ is sufficiently small. 
Therefore if $V'(x_0)$ is positive definite and $\Tr V'(x_0)\ge 2d\sigma^2$ we have 
\[
\int_{B_\delta(0)}\max\left\{\norm*{\pi'(x)},\norm*{\pi''(x)}\right\}\,\mu_\sigma(dx)\,<\,\infty, 
\] 
as desired. 
\end{proof}
\end{prop}
 
We conjecture that a similar proposition, with similar conditions on the drift and diffusion coefficients of \eqref{eq:SDE}, can be worked out for the hypotheses of Theorem \ref{thm:QuasiFrequency}. 
However, this would only guarantee that the convergence \eqref{eq:QuasiConvergence}~holds. 
It would not necessarily guarantee that $\left|\frac{1}{t}\pi_1(X_t^\sigma)-c_\sigma\right|$ would be small for values of $t$ which are less than $\tau_\sigma$ with high probability. 
That is, while an anologue of Proposition \ref{cor:}~may guarantee the existence of a quasi-asymptotic frequency, it would not guarantee its observability.

 

\end{document}